\documentclass[11pt]{article}%
\usepackage{amsfonts}
\usepackage{amsmath}
\usepackage{amssymb}
\usepackage{graphicx}%
\providecommand{\U}[1]{\protect\rule{.1in}{.1in}}
\everymath{\displaystyle}
\newtheorem{theorem}{Theorem}

\newtheorem{corollary}[theorem]{Corollary}

\newtheorem{definition}[theorem]{Definition}
\newtheorem{example}[theorem]{Example}

\newtheorem{lemma}[theorem]{Lemma}

\newtheorem{proposition}[theorem]{Proposition}

\newenvironment{proof}[1][Proof]{\noindent\textbf{#1.} }{\ \rule{0.5em}{0.5em}}
\begin{document}

\date{}
\title{Counting methods of area integrals and \textbf{Tchebychev }polynomials of
second kind on the ellipse}
\author{Abdelhamid Rehouma\\Department of mathematics Faculty of exact sciences\\University Hama Lakhdar, Eloued Algeria.\\URL: https://sites.google.com/view/mathsrehoumablog/accueil\\E-mail : rehoumaths@gmail.com}
\maketitle

\begin{abstract}
We use Gronwall's area formula to find the area of some differents regions as
circles, ellipses and lemniscates.We use Laurent and Taylor series expansions
of conformal mapping from the exterior of the unit disk to either of these
regions to compute the area of them.We close this work with the discussion of
orthogonal \textbf{Tchebychev }polynomials of second kind on the ellipse and interpolation.

\end{abstract}

\bigskip\textbf{Keywords:}{\scriptsize Unit disk,\ Conformal mapping ,Area
integral,Gronwall's area formula, univalents functions,regular functions
,Circles, Ellipses ,Lemniscates, Differentiation. \vspace{4mm} \newline}

\noindent\textbf{2000. MSC:}{\footnotesize \ 42C05, 33C45.}

\section{Introduction}

\bigskip In section two we use Gronwall's area formula to find the area of
some differents regions.We also give other methods to find of differents
examples of area integrals.

Suppose $G\subset%
\mathbb{C}
$ is an arbitrary domain ,$f$ \ is analytic in $G$ and set
\begin{equation}
I\left(  f\right)  =%
{\displaystyle\iint\limits_{G}}
\left\vert f\left(  z\right)  \right\vert ^{2}dm\text{ \ \ \ \ \ ,}z=x+iy
\label{Integral}%
\end{equation}

where $dm\left(  x,y\right)  =dxdy$ is the two-dimensional Lebesgue measure.

$I\left(  f\right)  $ can be represented explicitily in terms of the
coefficients of $f$ . Let us compute $I\left(  f\right)  $ for a special case
.Let $G=\left\{  z:r<\left\vert z\right\vert <R\right\}  ,0\leq r<R<\infty
,$and let
\[
f\left(  z\right)  =%
{\displaystyle\sum\limits_{n=-\infty}^{\infty}}
a_{n}z^{n}\text{ \ \ \ \ ,}z\in G
\]
it follows that,%
\[
I\left(  f\right)  =%
{\displaystyle\int\limits_{\rho=r}^{R}}
{\displaystyle\int\limits_{\phi=0}^{2\pi}}
\left(
{\displaystyle\sum\limits_{n=-\infty}^{\infty}}
a_{n}\rho^{n}e^{in\phi}\right)  \left(
{\displaystyle\sum\limits_{m=-\infty}^{\infty}}
\overline{a_{m}}\rho^{m}e^{-im\phi}\right)  \rho d\rho d\phi
\]%
\[
=2\pi%
{\displaystyle\int\limits_{\rho=r}^{R}}
{\displaystyle\sum\limits_{n=-\infty}^{\infty}}
\left\vert a_{n}\right\vert ^{2}\rho^{2n+1}d\rho
\]%
\[
=2\pi%
{\displaystyle\sum\limits_{n=-\infty}^{\infty}}
\left\vert a_{n}\right\vert ^{2}%
{\displaystyle\int\limits_{\rho=r}^{R}}
\rho^{2n+1}d\rho
\]
i-e%
\begin{equation}
I\left(  f\right)  =\pi%
{\displaystyle\sum\limits_{n=-\infty}^{\infty}}
\frac{R^{2n+2}-r^{2n+2}}{n+1}\left\vert a_{n}\right\vert ^{2} \label{Exp}%
\end{equation}
Since
\[
a_{0}=f\left(  0\right)
\]%
\[
a_{1}=f^{\prime}\left(  0\right)
\]%
\[
............
\]%
\[
a_{k}=\dfrac{f^{\left(  k\right)  }\left(  0\right)  }{k!}%
\]
We can state%
\begin{equation}
\left\vert f\left(  0\right)  \right\vert ^{2}\leq\frac{1}{R^{2}-r^{2}%
}I\left(  f\right)  \label{f1}%
\end{equation}
and%
\begin{equation}
\left\vert f^{\prime}\left(  0\right)  \right\vert ^{2}\leq\frac{2}%
{R^{4}-r^{4}}I\left(  f\right)  \label{f2}%
\end{equation}
Hence%
\begin{equation}
\left\vert f^{\left(  k\right)  }\left(  0\right)  \right\vert \leq
\frac{\left(  k+1\right)  !}{R^{2k+2}-r^{2k+2}}I\left(  f\right)  \label{fk}%
\end{equation}

\begin{definition}%
\begin{equation}
L^{2}\left(  G\right)  =\left\{  f:f\text{ analytic in }G\text{ and
\ }I\left(  f\right)  <\infty\right\}  \label{defi}%
\end{equation}
in this case,if $d_{z}=dist\left(  z,\partial G\right)  $%
\begin{equation}
\left\vert f\left(  z\right)  \right\vert ^{2}\leq\frac{I\left(  f\right)
}{\pi d_{z}^{2}} \label{inequ}%
\end{equation}

\end{definition}

\begin{theorem}
\textbf{Caratheodory} theorem \cite{Smirn},\cite{Smirno},\cite{Koosis}%
,\cite{Szeg},\cite{Gaier}.

The conformal mapping $\varphi$ of $G$ onto the disk $D_{R}=\left\{
w,\left\vert w\right\vert <R\right\}  $ such that%
\[
\varphi\left(  \xi\right)  =0,\text{ \ \ \ and \ \ }\varphi^{\prime}\left(
\xi\right)  =1
\]

can be conformally continued by \textbf{Caratheodory } theorem,\cite{Smirn}%
,\cite{Smirno},\cite{Koosis},\cite{Szeg},\cite{Gaier} If we

suppose that
\[
\psi:D_{R}=\left\{  w,\left\vert w\right\vert <R\right\}  \rightarrow G
\]

is its inverse map, then we have,
\[
\varphi\left(  t\right)  =u\Longleftrightarrow t=\psi\left(  u\right)  \text{
\ \ \ \ \ \ \ \ \ \ \ \ \ \ \ \ \ \ \ }\left(  t\in E\right)
\]

where
\[
\varphi\left(  t\right)  =\lim_{\underset{z\in G}{z\rightarrow t}}%
\varphi\left(  z\right)  \text{\ \ \ \ \ \ \ \ \ \ \ \ \ \ \ \ \ }\left(
\text{\ }t\in E\right)
\]

and%
\[
\psi\left(  e^{i\theta}\right)  =\lim_{r\longrightarrow1^{-}}\psi\left(
re^{i\theta}\right)  \ \ \ \ \ \ \ \ \ \ \left(  0\leq\theta\leq2\pi\right)
\]

These limit \ boundary values of $\varphi$ and $\psi$ are respectively
bijections and inversely from the contour $E$ onto unit circle $\left\{
u:\left\vert u\right\vert =1\right\}  .$
\end{theorem}

If $w=u+iv=f\left(  z\right)  $ $,z=x+iy\in D$, whose boundary is the closed
curve $C$,the Jacobien of the transformation , $u,v\longrightarrow x,y$ is
$\left\vert f^{\prime}\left(  z\right)  \right\vert ^{2}.$Hence the area of
$f\left(  D\right)  $ can be expressed as integral\textbf{\ :}%

\begin{equation}
A=%
{\displaystyle\iint\limits_{f\left(  D\right)  }}
dudv=%
{\displaystyle\iint\limits_{D}}
\left\vert f^{\prime}\left(  z\right)  \right\vert ^{2}dxdy \label{Areaa}%
\end{equation}
By Green's formula,\cite{Gaier},\cite{Zeev},\cite{Smirno}%
\begin{equation}%
{\displaystyle\iint\limits_{D}}
f\left(  z\right)  \overline{g^{\prime}\left(  z\right)  }dxdy=\frac{1}{2i}%
{\displaystyle\int\limits_{C}}
f\left(  z\right)  \overline{g\left(  z\right)  }dz \label{Greens}%
\end{equation}
Notice that%
\begin{equation}
A=%
{\displaystyle\iint\limits_{D}}
\left\vert f^{\prime}\left(  z\right)  \right\vert ^{2}dxdy=%
{\displaystyle\iint\limits_{D}}
f\left(  z\right)  \overline{f^{\prime}\left(  z\right)  }dxdy \label{result1}%
\end{equation}
ie%
\begin{equation}
A=%
{\displaystyle\iint\limits_{D}}
\left\vert f^{\prime}\left(  z\right)  \right\vert ^{2}dxdy=\frac{1}{2i}%
{\displaystyle\int\limits_{C}}
f^{\prime}\left(  z\right)  \overline{f\left(  z\right)  }dz \label{resultm}%
\end{equation}
If $f$ is regular and univalent in a domain $D$ and satisfies there
$\left\vert f\left(  z\right)  \right\vert \leq1$ then ,%
\begin{equation}%
{\displaystyle\iint\limits_{D}}
\left\vert f^{\prime}\left(  z\right)  \right\vert ^{2}dxdy\leq\pi
\label{formula}%
\end{equation}
Suppose that we are given the function analytique on the unit disk $D=\left\{
\left\vert z\right\vert \leq1\right\}  $ function $f\left(  z\right)  $ $=%
{\displaystyle\sum\limits_{n=0}^{\infty}}
a_{n}z^{n}$ .Regarding the magnitude of the coefficients can be obtanined by
means of Parseval's identity,,\cite{Gaier},\cite{Zeev},\cite{Smirno}%
\begin{equation}
\frac{1}{2\pi}%
{\displaystyle\int\limits_{0}^{2\pi}}
\left\vert f\left(  re^{i\theta}\right)  \right\vert ^{2}d\theta=%
{\displaystyle\sum\limits_{n=0}^{\infty}}
\left\vert a_{n}\right\vert ^{2}r^{2n} \label{Parseval}%
\end{equation}
which has also many other useful applications in the theory of functions.To
prove \eqref{Parseval},we observe that%
\[
\left\vert f\left(  re^{i\theta}\right)  \right\vert ^{2}=f\left(
re^{i\theta}\right)  \overline{f\left(  re^{i\theta}\right)  }=%
{\displaystyle\sum\limits_{n=0}^{\infty}}
a_{n}r^{n}e^{in\theta}%
{\displaystyle\sum\limits_{n=0}^{\infty}}
\overline{a_{n}}r^{n}e^{-in\theta}%
\]
i-e%
\[
\left\vert f\left(  re^{i\theta}\right)  \right\vert ^{2}=%
{\displaystyle\sum\limits_{n=0}^{\infty}}
{\displaystyle\sum\limits_{m=0}^{\infty}}
a_{n}\overline{a_{m}}r^{n+m}e^{i\left(  n-m\right)  \theta}%
\]
The rearrangement of the terms are permissible .Hence,%
\[%
{\displaystyle\int\limits_{0}^{2\pi}}
\left\vert f\left(  re^{i\theta}\right)  \right\vert ^{2}d\theta=%
{\displaystyle\sum\limits_{n=0}^{\infty}}
{\displaystyle\sum\limits_{m=0}^{\infty}}
a_{n}\overline{a_{m}}r^{n+m}%
{\displaystyle\int\limits_{0}^{2\pi}}
e^{i\left(  n-m\right)  \theta}d\theta
\]
reduces to
\begin{equation}%
{\displaystyle\int\limits_{0}^{2\pi}}
\left\vert f\left(  re^{i\theta}\right)  \right\vert ^{2}d\theta=2\pi%
{\displaystyle\sum\limits_{n=0}^{\infty}}
\left\vert a_{n}\right\vert ^{2}r^{2n} \label{identity}%
\end{equation}
If $\left\vert f\left(  z\right)  \right\vert \leq M\left(  r\right)  $ for
$\left\vert z\right\vert =r$, then obviously,,\cite{Gaier},\cite{Zeev}%
,\cite{Smirno}
\[%
{\displaystyle\sum\limits_{n=0}^{\infty}}
\left\vert a_{n}\right\vert ^{2}r^{2n}\leq M^{2}\left(  r\right)
\]
We known that \ If
\[
f\left(  z\right)  =%
{\displaystyle\sum\limits_{n=0}^{\infty}}
a_{n}z^{n}\text{ \ \ \ \ \ ,}\left\vert z\right\vert \leq r
\]
and%
\[
g\left(  z\right)  =%
{\displaystyle\sum\limits_{n=0}^{\infty}}
b_{n}z^{n}\text{ \ \ \ \ \ ,}\left\vert z\right\vert \leq\rho
\]
then,\cite{Gaier},\cite{Zeev},\cite{Smirno}
\begin{equation}
\frac{1}{2\pi i}%
{\displaystyle\int\limits_{\left\vert w\right\vert =r}}
f\left(  w\right)  g\left(  \frac{z}{w}\right)  \frac{dw}{w}=%
{\displaystyle\sum\limits_{n=0}^{\infty}}
a_{n}b_{n}z^{n}\text{ \ \ \ \ \ \ ,}\left\vert z\right\vert <\rho r
\label{identity2}%
\end{equation}
In polar coordinates $r,\theta$ of $\left\vert z\right\vert =r$ ,expressed as
$R,\phi$ in the $w$-plane ,when $f$ maps $\left\vert z\right\vert =r,$ $r<1$
in the positive direction , to the unit circle,$f\left(  re^{i\theta}\right)
=R\left(  r,\theta\right)  +i\phi\left(  r,\theta\right)  $
\begin{equation}
A=-\frac{1}{2}%
{\displaystyle\int\limits_{0}^{2\pi}}
R^{2}\left(  r,\theta\right)  \frac{\partial\phi}{\partial\theta}\left(
r,\theta\right)  d\theta\label{Polar}%
\end{equation}
Since by the Cauchy Riemann equations in polar coordinates $r,\theta$%
\begin{equation}
\frac{\partial\phi}{\partial\theta}\left(  r,\theta\right)  =\frac{r}{R}%
\frac{\partial R}{\partial r}\left(  r,\theta\right)  \label{Cauchy}%
\end{equation}
it follows that,\cite{Gaier},\cite{Zeev},\cite{Smirno}%
\begin{equation}
A=-\frac{r}{4}\frac{\partial}{\partial r}\left(
{\displaystyle\int\limits_{0}^{2\pi}}
\left\vert f\left(  re^{i\theta}\right)  \right\vert ^{2}d\theta\right)
\label{Arear}%
\end{equation}

\begin{example}
Let $M$ is a region whose boundary \ is the cardioid has a parametric polar
equation,for $0\leq\theta\leq2\pi$%
\[
r=\frac{1}{2}\left(  1+\cos\theta\right)
\]

The area $A_{cardioid}$ can be found by integration .We have%
\[
A_{cardioid}=%
{\displaystyle\int\limits_{0}^{\pi}}
r^{2}d\theta=\frac{1}{4}%
{\displaystyle\int\limits_{0}^{\pi}}
\left(  1+\cos\theta\right)  ^{2}d\theta
\]%
\[
=\frac{1}{4}%
{\displaystyle\int\limits_{0}^{\pi}}
\left(  \frac{3}{2}+2\cos\theta+\frac{1}{2}\cos2\theta\right)  d\theta
=\frac{3}{8}\pi
\]

\end{example}

\begin{proposition}
The class of univalents functions in $\left\vert z\right\vert <1$ which have a
Laurent expansion,\cite{Gaier},\cite{Zeev},\cite{Smirno}%
\begin{equation}
f\left(  z\right)  =\frac{1}{z}+%
{\displaystyle\sum\limits_{n=0}^{\infty}}
b_{n}z^{n} \label{Laurent}%
\end{equation}

The coefficients $b_{0},b_{1},b_{2}.....b_{n}....$are subject to the
inequality
\begin{equation}%
{\displaystyle\sum\limits_{n=0}^{\infty}}
n\left\vert b_{n}\right\vert ^{2}\leq1 \label{Inequality}%
\end{equation}

\end{proposition}

\begin{proof}%
\[%
{\displaystyle\int\limits_{0}^{2\pi}}
\left\vert f\left(  re^{i\theta}\right)  \right\vert ^{2}d\theta=%
{\displaystyle\int\limits_{0}^{2\pi}}
f\left(  re^{i\theta}\right)  \overline{f\left(  re^{i\theta}\right)  }d\theta
\]%
\[
=%
{\displaystyle\int\limits_{0}^{2\pi}}
\left[  \frac{1}{re^{i\theta}}+%
{\displaystyle\sum\limits_{n=0}^{\infty}}
b_{n}r^{n}e^{in\theta}\right]  \left[  \frac{1}{re^{-i\theta}}+%
{\displaystyle\sum\limits_{n=0}^{\infty}}
b_{n}r^{n}e^{-in\theta}\right]  d\theta
\]%
\[
=2\pi\left[  \frac{1}{r^{2}}+%
{\displaystyle\sum\limits_{n=0}^{\infty}}
\left\vert b_{n}\right\vert ^{2}r^{2n}\right]
\]

Hence,%
\[
A=-\frac{r}{4}\frac{\partial}{\partial r}\left(
{\displaystyle\int\limits_{0}^{2\pi}}
\left\vert f\left(  re^{i\theta}\right)  \right\vert ^{2}d\theta\right)
\]

then%

\[
A=-\frac{\pi r}{2}\frac{\partial}{\partial r}\left[  \frac{1}{r^{2}}+%
{\displaystyle\sum\limits_{n=0}^{\infty}}
\left\vert b_{n}\right\vert ^{2}r^{2n}\right]
\]

and thus%
\[
\frac{1}{\pi}A=\frac{1}{r^{2}}-%
{\displaystyle\sum\limits_{n=1}^{\infty}}
n\left\vert b_{n}\right\vert ^{2}r^{2n}%
\]

Since the area can not negative it follows that
\[%
{\displaystyle\sum\limits_{n=1}^{\infty}}
n\left\vert b_{n}\right\vert ^{2}r^{2n}\leq\frac{1}{r^{2}}%
\]

holds for all values of $r$ between $0$ and $1,$see,\cite{Gaier}%
,\cite{Zeev},\cite{Smirno}.It therefore must also be true for
$r\longrightarrow1.$This proves ,\eqref{Inequality}.
\end{proof}

\section{Gronwall's area formula}

\begin{lemma}
Suppose that a simple closed curve $C$ is defined parametrically in the
complex w-plane
\begin{equation}
u=u\left(  \theta\right)  \text{ \ \ ,}v=v\left(  \theta\right)  \text{
\ \ \ ,}0\leq\theta\leq2\pi\label{Param}%
\end{equation}

Then the area $A$ enclosed by $C$ is given by,\cite{Gaier},\cite{Zeev}%
,\cite{Gradsh}.
\begin{equation}
A=\frac{1}{2}%
{\displaystyle\int\limits_{0}^{2\pi}}
\left(  u\frac{dv}{d\theta}-v\frac{du}{d\theta}\right)  d\theta\label{Gronw}%
\end{equation}

\end{lemma}

\begin{proof}
By Green's theorem ,if $p=p\left(  u,v\right)  $ and $q=q\left(  u,v\right)  $
are functions with continuous partial derivatives on a region $B$ enclosed by
$C$ ,\cite{Gaier},\cite{Zeev},\cite{Smirno}then%

\[
A=%
{\displaystyle\iint\limits_{B}}
\left(  \frac{\partial p}{\partial u}+\frac{\partial q}{\partial v}\right)
dudv=%
{\displaystyle\int\limits_{C}}
\left(  -qdu+pdv\right)
\]

Choosing ,$p=u$ ,$q=v$ gives%

\[
A=%
{\displaystyle\iint\limits_{B}}
dudv=\frac{1}{2}%
{\displaystyle\int\limits_{C}}
\left(  -vdu+udv\right)  =\frac{1}{2}%
{\displaystyle\int\limits_{0}^{2\pi}}
\left(  u\frac{dv}{d\theta}-v\frac{du}{d\theta}\right)  d\theta
\]

We shall use this result to prove Gronwall's area theorem.
\end{proof}

\begin{theorem}
If there is a conformal mapping $\Psi_{r}$ ,from the exterior of $D_{r}$
,$r>1$ to the exterior of $C_{r}$ of the form%

\begin{equation}
\Psi_{r}\left(  z\right)  =z+%
{\displaystyle\sum\limits_{n=0}^{\infty}}
\frac{b_{n}}{z^{n}} \label{conform}%
\end{equation}

then the area $A_{r\text{ }}$of \ the region $B_{r}$ enclosed by $C_{r}$ is
given by,\cite{Gaier},\cite{Zeev},\cite{Smirno}%
\begin{equation}
A_{r\text{ }}=\pi\left(  r^{2}-%
{\displaystyle\sum\limits_{n=1}^{\infty}}
n\left\vert b_{n}\right\vert ^{2}r^{-2n}\right)  \label{ConfArea}%
\end{equation}

\end{theorem}

\begin{proof}
Firstly,we write the mapping $\Psi_{r}$ \ as%
\[
w=\Psi_{r}\left(  z\right)  =u\left(  r,\theta\right)  +iv\left(
r,\theta\right)  \text{ \ \ \ \ \ \ \ \ , \ \ }z=re^{i\theta}%
\]

where $u$ and $v$ are real-valued functions.By \eqref{Gronw} and since $r$ is
constant, the area enclosed by $C_{r}$ is given by%
\begin{equation}
A_{r}=\frac{1}{2}%
{\displaystyle\int\limits_{0}^{2\pi}}
\left(  u\frac{dv}{d\theta}-v\frac{du}{d\theta}\right)  d\theta\label{Ar1}%
\end{equation}

We now want to prove that%
\begin{equation}
A_{r}=\operatorname{Im}\left\{  \frac{1}{2}%
{\displaystyle\int\limits_{0}^{2\pi}}
\overline{\Psi_{r}\left(  z\right)  }\Psi_{r}^{\prime}\left(  z\right)
dz\right\}  \label{Ar2}%
\end{equation}

where $\operatorname{Im}$ denote the imaginary part .Since%
\begin{equation}
\Psi_{r}^{\prime}\left(  z\right)  =e^{-i\theta}\left(  \frac{\partial
u}{\partial r}+i\frac{\partial v}{\partial r}\right)  =\frac{e^{-i\theta}}%
{r}\left(  \frac{\partial v}{\partial\theta}-i\frac{\partial u}{\partial
\theta}\right)  \label{Derivpsci}%
\end{equation}

on using the Cauchy Riemann equations in polar coordinates form,.In fact
$dz=ire^{i\theta}d\theta,$we see that%
\[
\operatorname{Im}\left\{  \frac{1}{2}%
{\displaystyle\int\limits_{0}^{2\pi}}
\overline{\Psi_{r}\left(  z\right)  }\Psi_{r}^{\prime}\left(  z\right)
dz\right\}  =\operatorname{Im}\left\{  \frac{1}{2}%
{\displaystyle\int\limits_{0}^{2\pi}}
\left(  u-iv\right)  \left(  \frac{\partial v}{\partial\theta}-i\frac{\partial
u}{\partial\theta}\right)  d\theta\right\}
\]%
\[
=\frac{1}{2}%
{\displaystyle\int\limits_{0}^{2\pi}}
\left(  u\frac{\partial v}{\partial\theta}-v\frac{\partial u}{\partial\theta
}\right)  d\theta=A_{r}%
\]

We can now rewrite equation \eqref{conform} as
\begin{equation}
\Psi_{r}\left(  z\right)  =%
{\displaystyle\sum\limits_{n=-1}^{\infty}}
\frac{b_{n}}{z^{n}}\text{ \ \ \ \ \ , }b_{-1}=1 \label{equation}%
\end{equation}

Using \eqref{equation} and \eqref{Ar2} we have%
\[
A_{r}=\operatorname{Im}\left\{  \frac{1}{2}%
{\displaystyle\int\limits_{0}^{2\pi}}
\left(
{\displaystyle\sum\limits_{n=-1}^{\infty}}
\frac{\overline{b}_{n}}{\overline{z}^{n}}%
{\displaystyle\sum\limits_{m=-1}^{\infty}}
\frac{-mb_{m}}{z^{m+1}}\right)  dz\right\}
\]%
\[
=\operatorname{Im}\left\{  -\frac{1}{2}%
{\displaystyle\sum\limits_{m=-1}^{\infty}}
{\displaystyle\sum\limits_{n=-1}^{\infty}}
mb_{m}\overline{b}_{n}%
{\displaystyle\int\limits_{\left\vert z=r\right\vert }}
\overline{z}^{-n}z^{m-1}dz\right\}
\]

By writing $z=re^{i\theta}$,we find%
\[%
{\displaystyle\int\limits_{\left\vert z=r\right\vert }}
\overline{z}^{-n}z^{m-1}dz=ir^{-m-n}%
{\displaystyle\int\limits_{0}^{2\pi}}
e^{i\left(  n-m\right)  }d\theta=0\text{ \ \ \ , }n\neq m
\]

and%
\[%
{\displaystyle\int\limits_{\left\vert z=r\right\vert }}
\overline{z}^{-n}z^{n-1}dz=2\pi ir^{-2n}%
\]

therefore%
\[
A_{r}=\operatorname{Im}\left\{  -i\pi%
{\displaystyle\sum\limits_{n=-1}^{\infty}}
n\left\vert b_{n}\right\vert ^{2}r^{-2n}\right\}  =-\pi%
{\displaystyle\sum\limits_{n=-1}^{\infty}}
n\left\vert b_{n}\right\vert ^{2}r^{-2n}%
\]%
\[
=\pi\left(  r^{2}-%
{\displaystyle\sum\limits_{n=1}^{\infty}}
n\left\vert b_{n}\right\vert ^{2}r^{-2n}\right)
\]

The last relation proves the theorem.
\end{proof}

\section{\bigskip Circles, ellipses and lemniscates}

\bigskip In this case $C_{r}$ is the circle with radius $r$ centred ot the
origin lying in the complex $w$-plane .The conformal mapping from the exterior
of $D_{r}$ onto the exterior of $C_{r}$ is given by :%
\begin{equation}
w=\Psi_{r}\left(  z\right)  =z \label{confor1}%
\end{equation}

Thus we see that all the $b_{n}$'s from equation\eqref{conform} are zero.By
Gronwall's area theorem,the area $A_{r}$ of the circle $C_{r}$ is given by
\begin{equation}
A_{r}=\pi\left(  r^{2}-%
{\displaystyle\sum\limits_{n=1}^{\infty}}
n\left\vert b_{n}\right\vert ^{2}r^{-2n}\right)  =\pi r^{2} \label{Circle}%
\end{equation}

\bigskip The conformal mapping from the exterior of the unit circle centred at
the origin of the complex $z-$plane $D_{r}$ onto an ellipse $C_{r}$ in the
complex $w$-plane with semi axis $a=r+\frac{1}{r}$ ,$b=r-\frac{1}{r}$ is

given by%
\begin{equation}
w=\Psi_{r}\left(  z\right)  =z+\frac{1}{z} \label{coonfor2}%
\end{equation}

Thus we see that all the $b_{n}$'s from equation\eqref{conform} are zero
except for $b_{1}=1.$By Gronwall's area theorem,the area $A_{r}$ of the circle
$C_{r}$ is given by%
\[
A_{r}=\pi\left(  r^{2}-%
{\displaystyle\sum\limits_{n=1}^{\infty}}
n\left\vert b_{n}\right\vert ^{2}r^{-2n}\right)  =\pi\left(  r+\frac{1}%
{r}\right)  \left(  r-\frac{1}{r}\right)
\]

i-e%
\begin{equation}
A_{r}=\pi ab \label{Ellipse}%
\end{equation}

which agrees with the classical formula for the area of an ellipse.

A lemniscate is shapped with any number of leaves .The boundary of $m$-leafed
symmetric lemniscate is the set%
\begin{equation}
\left\{  w\in%
\mathbb{C}
,\left\vert w^{m}-1\right\vert =0\right\}  ,m=2,3,4..... \label{lemnisc}%
\end{equation}

The conformal mapping from the exterior of the uniit circle centred at the
origin of the complex $z-$plane $D_{r}$ onto the exterior of the $m$-leafed
symmetric lemniscate\ in the complex $w-$plane is

given by
\begin{equation}
w=\Psi\left(  z\right)  =z\left(  1+\frac{1}{z^{m}}\right)  ^{\frac{1}{m}}
\label{confor3}%
\end{equation}

The image of this mapping to include the boundary of the lemniscate ,we look
at the image on the unit circle in the the complex $z-$plane, that is, take
$r=1.$In this case
\[
w^{m}=z^{m}+1
\]

so that
\[
\left\vert w^{m}-1\right\vert =\left\vert z\right\vert ^{m}=1
\]

\begin{theorem}
If $z$ describe the unit circle.Then the area $A$ of $\ $an $m$-leafed
symetric lemniscate is given by%
\begin{equation}
A=\pi\left(
{\displaystyle\sum\limits_{n=0}^{\infty}}
\left(  C_{n}^{\frac{1}{m}}\right)  ^{2}-m%
{\displaystyle\sum\limits_{n=1}^{\infty}}
n\left(  C_{n}^{\frac{1}{m}}\right)  ^{2}\right)  \label{lemniscate}%
\end{equation}

Hence%
\begin{equation}
A=2^{\frac{2}{m}-1}\frac{\Gamma\left(  \frac{1}{m}+\frac{1}{2}\right)
}{\Gamma\left(  \frac{1}{m}+1\right)  }\sqrt{\pi} \label{lemniscate1}%
\end{equation}

\end{theorem}

\begin{proof}
If $z$ describe the unit circle.Then by \eqref{confor3} for large $\left\vert
z\right\vert $%
\begin{equation}
\Psi\left(  z\right)  =z%
{\displaystyle\sum\limits_{n=0}^{\infty}}
C_{n}^{\frac{1}{m}}\frac{1}{z^{mn-1}} \label{confor5}%
\end{equation}

According to Gronwall's area formula,\eqref{ConfArea} the area $A$ of $\ $an
$m$-leafed symetric lemniscate is given by%
\[
A=\pi\left[  1-%
{\displaystyle\sum\limits_{n=1}^{\infty}}
n\left\vert b_{n}\right\vert ^{2}\right]
\]%
\[
=\pi\left[  1-%
{\displaystyle\sum\limits_{n=1}^{\infty}}
\left(  mn-1\right)  \left\vert b_{mn-1}\right\vert ^{2}\right]
\]

by \eqref{confor5}
\[
A=\pi m%
{\displaystyle\sum\limits_{n=1}^{\infty}}
\left(  \frac{1}{m}-n\right)  \left(  C_{n}^{\frac{1}{m}}\right)  ^{2}%
\]%
\[
=\pi\left(
{\displaystyle\sum\limits_{n=0}^{\infty}}
\left(  C_{n}^{\frac{1}{m}}\right)  ^{2}-m%
{\displaystyle\sum\limits_{n=1}^{\infty}}
n\left(  C_{n}^{\frac{1}{m}}\right)  ^{2}\right)
\]

and \eqref{lemniscate} is proved.To prove \eqref{lemniscate1}, we known ( see
Gradshteyn and Ryzhik,\cite{Gradsh} page 5 ) ,we known that%
\begin{equation}%
{\displaystyle\sum\limits_{k=0}^{\infty}}
\left(  C_{k}^{\alpha}\right)  ^{2}=\frac{\Gamma\left(  2\alpha+1\right)
}{\left(  \Gamma\left(  \alpha+1\right)  \right)  ^{2}} \label{Ryzhi1}%
\end{equation}
and%
\begin{equation}%
{\displaystyle\sum\limits_{k=0}^{\infty}}
k\left(  C_{k}^{\alpha}\right)  ^{2}=\frac{\Gamma\left(  2\alpha\right)
}{\left(  \Gamma\left(  \alpha\right)  \right)  ^{2}} \label{Ryzhi2}%
\end{equation}

Using these assumptions and proprietes of the Gamma function we find
\[
A=\pi\left(  \frac{\Gamma\left(  \frac{2}{m}+1\right)  }{\left(  \Gamma\left(
\frac{1}{m}+1\right)  \right)  ^{2}}-m\frac{\Gamma\left(  \frac{2}{m}\right)
}{\left(  \Gamma\left(  \frac{1}{m}\right)  \right)  ^{2}}\right)
\]

By the Duplicate formula ( see Gradshteyn and Ryzhik, \cite{Gradsh} page 946
),we get
\[
A=2^{\frac{2}{m}-1}\frac{\Gamma\left(  \frac{1}{m}+\frac{1}{2}\right)
}{\Gamma\left(  \frac{1}{m}+1\right)  }\sqrt{\pi}%
\]

We will now confirm this result by using polar coordinates. Writing $w=\rho
e^{i\phi},$we have
\[
\left\vert w^{m}-1\right\vert =1
\]
or%
\[
\left(  w^{m}-1\right)  \left(  \overline{w}^{m}-1\right)
\]
i-e%
\[
\rho^{2m}-2\rho^{m}\cos m\phi=0
\]

so that
\[
\rho=2^{\frac{1}{m}}\cos^{\frac{1}{m}}\left(  m\phi\right)
\]

The area of one leaf of the $m$-leafed lemniscate is given by,\cite{Gradsh}
\[
\frac{1}{m}A=\frac{1}{2}%
{\displaystyle\int\limits_{-\frac{\pi}{2m}}^{\frac{\pi}{2m}}}
\rho^{2}d\phi
\]

and thus the area of the $m$-leafed lemniscate is given by,\cite{Gradsh}%
\[
A=m%
{\displaystyle\int\limits_{0}^{\frac{\pi}{2m}}}
\rho^{2}d\phi
\]%
\[
=2^{\frac{2}{m}}m%
{\displaystyle\int\limits_{0}^{\frac{\pi}{2m}}}
\cos^{\frac{2}{m}}\left(  m\phi\right)  d\phi
\]

If $u=\cos m\phi$,then%
\[
A=2^{\frac{2}{m}}%
{\displaystyle\int\limits_{0}^{1}}
\frac{u^{\frac{2}{m}}}{\sqrt{1-u^{2}}}du
\]

If $v=u^{2},$then%
\[
A=2^{\frac{2}{m}-1}%
{\displaystyle\int\limits_{0}^{1}}
v^{\frac{1}{m}-\frac{1}{2}}\left(  1-v\right)  ^{-\frac{1}{2}}dv
\]

and therefore,by the Beta function,%
\[
A=2^{\frac{2}{m}-1}\frac{\Gamma\left(  \frac{1}{m}+\frac{1}{2}\right)
\Gamma\left(  \frac{1}{2}\right)  }{\Gamma\left(  \frac{1}{m}+\frac{1}%
{2}+\frac{1}{2}\right)  }%
\]%
\[
=2^{\frac{2}{m}-1}\frac{\Gamma\left(  \frac{1}{m}+\frac{1}{2}\right)  }%
{\Gamma\left(  \frac{1}{m}+1\right)  }\sqrt{\pi}%
\]

and the theorem is proved.
\end{proof}

In some cases to be considered ,we will take the limit $r\overset
{>}{\longrightarrow}1$ from above ,and assume that the image of the unit
circle is the boundary of the region for which we are trying to find the area.

\begin{proposition}
Consider the area integral\textbf{\ }%
\begin{equation}
\Im\left(  f\right)  =%
{\displaystyle\int\limits_{0}^{2\pi}}
d\varphi%
{\displaystyle\int\limits_{0}^{2\pi}}
\frac{\left\vert f\left(  e^{i\varphi}\right)  -f\left(  e^{i\theta}\right)
\right\vert ^{2}}{\left\vert e^{i\varphi}-e^{i\theta}\right\vert ^{2}}%
d\theta\label{Areaaa}%
\end{equation}
If%
\begin{equation}
\Im_{r}\left(  f\right)  =%
{\displaystyle\int\limits_{C_{r}}}
dy%
{\displaystyle\int\limits_{C_{1}}}
\frac{\left\vert f\left(  y\right)  -f\left(  x\right)  \right\vert ^{2}%
}{\left\vert y-x\right\vert ^{2}}dx \label{area 2}%
\end{equation}

where ,$C_{r}=\left\{  \left\vert z\right\vert =r,0\leq r\leq1\right\}
,$oriented in the positive direction then%
\begin{equation}
Lim_{r\overset{<}{\longrightarrow}1}\Im_{r}\left(  f\right)  =\Im\left(
f\right)  \label{limit1}%
\end{equation}

if $r<1.$Then,%
\begin{equation}%
{\displaystyle\int\limits_{C_{r}}}
dy%
{\displaystyle\int\limits_{C_{1}}}
\left\vert f\left(  y\right)  \right\vert ^{2}\frac{dx}{\left(  y-x\right)
^{2}}=0\text{ \ \ \ \ \ \ ,}r<1 \label{form1}%
\end{equation}%
\begin{equation}%
{\displaystyle\int\limits_{C_{r}}}
dy%
{\displaystyle\int\limits_{C_{1}}}
\left\vert f\left(  x\right)  \right\vert ^{2}\frac{dx}{\left(  y-x\right)
^{2}}=0\text{ } \label{form2}%
\end{equation}%
\begin{equation}%
{\displaystyle\int\limits_{C_{r}}}
dy\text{ }%
{\displaystyle\int\limits_{C_{r}}}
\overline{f\left(  x\right)  }f\left(  y\right)  \frac{dx}{\left(  y-x\right)
^{2}}=0\text{ \ \ \ \ }\left(  \left\vert z\right\vert =1,\left\vert
w\right\vert <1\right)  \label{form3}%
\end{equation}
If%
\begin{equation}
K_{r}\left(  f\right)  =%
{\displaystyle\int\limits_{C_{r}}}
dy%
{\displaystyle\int\limits_{C_{1}}}
\overline{f\left(  y\right)  }f\left(  x\right)  \frac{dx}{\left(  y-x\right)
^{2}} \label{krf}%
\end{equation}

then%
\begin{equation}
\Im\left(  f\right)  =4\pi^{2}%
{\displaystyle\sum\limits_{p=0}^{\infty}}
p\left\vert a_{p}\right\vert ^{2} \label{JrKrf}%
\end{equation}

and%
\begin{equation}
A=\frac{1}{4\pi}\Im\left(  f\right)  \label{form4}%
\end{equation}

\end{proposition}

\begin{proof}
The function $f$ \ is holomorphe on the disk $\left\vert z\right\vert \leq
1,$then the function%
\[
F\left(  \theta,\varphi\right)  =\frac{\left\vert f\left(  e^{i\varphi
}\right)  -f\left(  e^{i\theta}\right)  \right\vert ^{2}}{\left\vert
e^{i\varphi}-e^{i\theta}\right\vert ^{2}}%
\]

is continued on $\left[  0\text{ \ \ }2\pi\right]  \times\left[  0\text{
\ \ }2\pi\right]  $,then $\Im\left(  f\right)  $ exist.Setting,$x=e^{i\theta}%
$,$y=e^{i\varphi}$ we have%
\[
\left\vert e^{i\varphi}-e^{i\theta}\right\vert ^{2}=\left(  e^{i\varphi
}-e^{i\theta}\right)  \left(  e^{-i\varphi}-e^{-i\theta}\right)
=-e^{-i\left(  \theta+\varphi\right)  }\left(  y-x\right)  ^{2}%
\]%
\[
dx=ie^{i\theta}d\theta,dy=ie^{i\varphi}d\varphi
\]

then%
\[%
{\displaystyle\int\limits_{0}^{2\pi}}
d\varphi%
{\displaystyle\int\limits_{0}^{2\pi}}
\frac{\left\vert f\left(  e^{i\varphi}\right)  -f\left(  e^{i\theta}\right)
\right\vert ^{2}}{\left\vert e^{i\varphi}-e^{i\theta}\right\vert ^{2}}d\theta=%
{\displaystyle\int\limits_{C_{1}}}
dy%
{\displaystyle\int\limits_{C_{1}}}
\frac{\left\vert f\left(  y\right)  -f\left(  x\right)  \right\vert ^{2}%
}{\left\vert y-x\right\vert ^{2}}dx
\]

If%
\[
\Im_{r}\left(  f\right)  =-%
{\displaystyle\int\limits_{0}^{2\pi}}
d\varphi%
{\displaystyle\int\limits_{0}^{2\pi}}
\frac{\left\vert f\left(  re^{i\varphi}\right)  -f\left(  e^{i\theta}\right)
\right\vert ^{2}}{\left\vert re^{i\varphi}-e^{i\theta}\right\vert ^{2}%
}re^{i\left(  \theta+\varphi\right)  }d\theta
\]

then%
\[
Lim_{r\overset{<}{\longrightarrow}1}%
{\displaystyle\int\limits_{C_{r}}}
dy%
{\displaystyle\int\limits_{C_{1}}}
\frac{\left\vert f\left(  y\right)  -f\left(  x\right)  \right\vert ^{2}%
}{\left\vert y-x\right\vert ^{2}}dx=%
{\displaystyle\int\limits_{C_{1}}}
dy%
{\displaystyle\int\limits_{C_{1}}}
\frac{\left\vert f\left(  y\right)  -f\left(  x\right)  \right\vert ^{2}%
}{\left\vert y-x\right\vert ^{2}}dx
\]

i-e%
\[
Lim_{r\overset{<}{\longrightarrow}1}\Im_{r}\left(  f\right)  =\Im\left(
f\right)
\]

we have
\[%
{\displaystyle\int\limits_{C_{1}}}
\frac{dx}{\left(  y-x\right)  ^{2}}=0\text{ \ \ \ \ ,}\left\vert y\right\vert
<1
\]

and%
\[%
{\displaystyle\int\limits_{C_{r}}}
\frac{dy}{\left(  y-x\right)  ^{2}}=0\text{ \ \ \ \ ,}\left\vert x\right\vert
=1
\]

also we have%
\[%
{\displaystyle\int\limits_{C_{r}}}
f\left(  y\right)  \frac{dy}{\left(  y-x\right)  ^{2}}=0\text{ \ \ \ \ ,}%
\left\vert x\right\vert =1
\]

thus proprietes \eqref{form1},\eqref{form2},\eqref{form3} are proved.We have,%

\[
\left\vert f\left(  y\right)  -f\left(  x\right)  \right\vert ^{2}=\left\vert
f\left(  y\right)  \right\vert ^{2}+\left\vert f\left(  x\right)  \right\vert
^{2}-\overline{f\left(  y\right)  }f\left(  x\right)  -\overline{f\left(
x\right)  }f\left(  y\right)
\]

then,in view \eqref{form1},\eqref{form2},\eqref{form3} we get%
\[
\Im_{r}\left(  f\right)  =-%
{\displaystyle\int\limits_{C_{r}}}
dy%
{\displaystyle\int\limits_{C_{1}}}
\overline{f\left(  y\right)  }f\left(  x\right)  \frac{dx}{\left(  y-x\right)
^{2}}=-K_{r}\left(  f\right)
\]

By residu's theorem at the pole $x=y,$ we have%
\[%
{\displaystyle\int\limits_{C_{1}}}
f\left(  x\right)  \frac{dx}{\left(  y-x\right)  ^{2}}=2i\pi f^{\prime}\left(
y\right)
\]

then if $\left\vert y\right\vert =r<1,$we get%
\begin{equation}
\Im_{r}\left(  f\right)  =-K_{r}\left(  f\right)  =-2i\pi%
{\displaystyle\int\limits_{C_{r}}}
\overline{f\left(  y\right)  }f^{\prime}\left(  y\right)  dy \label{result}%
\end{equation}%
\[
=-2i\pi%
{\displaystyle\int\limits_{C_{r}}}
\overline{%
{\displaystyle\sum\limits_{q=0}^{\infty}}
r^{q}e^{i\varphi}}%
{\displaystyle\sum\limits_{p=0}^{\infty}}
pa_{p}r^{p}e^{i\theta}d\theta
\]%
\[
=2\pi%
{\displaystyle\int\limits_{0}^{2\pi}}
{\displaystyle\sum\limits_{p=0}^{\infty}}
{\displaystyle\sum\limits_{q=0}^{\infty}}
pa_{p}\overline{a_{q}}r^{p+q}e^{i\left(  p-q\right)  \varphi}d\varphi
\]%
\[
=4\pi^{2}%
{\displaystyle\sum\limits_{p=0}^{\infty}}
p\left\vert a_{p}\right\vert ^{2}r^{2p}%
\]

Thus%
\[
Lim_{r\overset{<}{\longrightarrow}1}\Im_{r}\left(  f\right)  =4\pi^{2}%
{\displaystyle\sum\limits_{p=0}^{\infty}}
p\left\vert a_{p}\right\vert ^{2}%
\]

By \eqref{result},and \eqref{Greens},we get
\[
Lim_{r\overset{<}{\longrightarrow}1}\Im_{r}\left(  f\right)  =\frac{4\pi}{2i}%
{\displaystyle\int\limits_{C_{1}}}
f^{\prime}\left(  y\right)  \overline{f\left(  y\right)  }dy=4\pi A=\Im\left(
f\right)
\]

i-e%
\[
A=\frac{1}{4\pi}\Im\left(  f\right)
\]

The proposition is proved.
\end{proof}

\begin{corollary}
Let $f\left(  z\right)  =%
{\displaystyle\sum\limits_{n=1}^{\infty}}
a_{n}z^{n},$on the disk $\left\vert z\right\vert \leq1,$if $b_{1}%
,b_{2}.....b_{n}....$are the coefficients of the power series expansion
\begin{equation}
z\frac{f^{\prime}\left(  z\right)  }{f\left(  z\right)  }=1+b_{1}z+b_{2}%
z^{2}+.... \label{condition}%
\end{equation}

Then%
\begin{equation}%
{\displaystyle\iint\limits_{D}}
\left\vert f^{\prime}\left(  z\right)  +zf^{\prime\prime}\left(  z\right)
\right\vert ^{2}dxdy=%
{\displaystyle\sum\limits_{p=1}^{\infty}}
p\left\vert a_{p}+b_{1}a_{p-1}+......b_{p-1}\right\vert ^{2} \label{resultt}%
\end{equation}

\end{corollary}

\begin{proof}
We have now,if $f\left(  z\right)  =%
{\displaystyle\sum\limits_{n=1}^{\infty}}
a_{n}z^{n}$%
\[
zf^{\prime}\left(  z\right)  =%
{\displaystyle\sum\limits_{n=1}^{\infty}}
na_{n}z^{n}%
\]

and%
\[
z\frac{f^{\prime}\left(  z\right)  }{f\left(  z\right)  }.f\left(  z\right)
=\left(  1+%
{\displaystyle\sum\limits_{n=1}^{\infty}}
b_{n}z^{n}\right)
{\displaystyle\sum\limits_{n=1}^{\infty}}
a_{n}z^{n}=%
{\displaystyle\sum\limits_{n=1}^{\infty}}
c_{n}z^{n}%
\]

where%
\[
c_{n}=a_{n}+b_{1}a_{n-1}+......b_{n-1}%
\]

Hence%
\[
\Im\left(  zf^{\prime}\left(  z\right)  \right)  =4\pi^{2}%
{\displaystyle\sum\limits_{p=1}^{\infty}}
p\left\vert a_{p}+b_{1}a_{p-1}+......b_{p-1}\right\vert ^{2}%
\]

This proves ,\eqref{resultt}.
\end{proof}

\begin{proposition}
Let $z_{1},z_{2},z_{3}....z_{p}$ be $n$ distincts points which are situated in
the interior of $D=\left\{  \left\vert z\right\vert <1\right\}  ,$
\begin{equation}
I_{p}=\frac{1}{\pi}%
{\displaystyle\iint\limits_{D}}
\left\vert \frac{1}{z-z_{p}}\right\vert ^{2}dxdy\text{\ \ \ \ , \ \ \ }z=x+iy
\label{Ip}%
\end{equation}

and%
\begin{equation}
J_{p}=\frac{1}{2\pi}%
{\displaystyle\iint\limits_{D}}
\left\vert \frac{z+z_{p}}{z-z_{p}}\right\vert ^{2}dxdy\text{ \ \ \ \ ,
\ \ \ }z=x+iy \label{Jp}%
\end{equation}

then%
\begin{equation}
I_{p}=\log\left\vert \frac{1-\left\vert z_{p}\right\vert ^{2}}{\left\vert
z_{p}\right\vert ^{2}}\right\vert \label{Ipvalue}%
\end{equation}

and%
\begin{equation}
J_{p}=\frac{1}{2}+2\left\vert z_{p}\right\vert ^{2}\log\left\vert
\frac{1-\left\vert z_{p}\right\vert ^{2}}{\left\vert z_{p}\right\vert ^{2}%
}\right\vert \label{Jpvalue}%
\end{equation}

therefore%
\begin{equation}%
{\displaystyle\iint\limits_{D}}
\left(
{\displaystyle\sum\limits_{p=1}^{m}}
\left\vert \frac{1}{z-z_{p}}\right\vert ^{2}\right)  dxdy=\pi\log\left\vert
{\displaystyle\prod\limits_{p=1}^{m}}
\frac{1-\left\vert z_{p}\right\vert ^{2}}{\left\vert z_{p}\right\vert ^{2}%
}\right\vert \label{SommeIp}%
\end{equation}

and%
\begin{equation}%
{\displaystyle\iint\limits_{D}}
\left(
{\displaystyle\sum\limits_{p=1}^{m}}
\left\vert \frac{z+z_{p}}{z-z_{p}}\right\vert ^{2}\right)  dxdy=m\pi+2\pi
\log\left\vert
{\displaystyle\prod\limits_{p=1}^{m}}
\frac{1-\left\vert z_{p}\right\vert ^{2}}{\left\vert z_{p}\right\vert ^{2}%
}\right\vert ^{2\left\vert z_{p}\right\vert ^{2}} \label{SommeJp}%
\end{equation}

\end{proposition}

\begin{proof}
suppose \ $z=re^{i\theta}=x+iy$%
\[
dxdy=rdrd\theta,\text{ \ }0\leq r<1,\text{ }0\leq\theta\leq2\pi
\]

then%
\[
I_{p}=\frac{1}{\pi}%
{\displaystyle\int\limits_{0}^{1}}
{\displaystyle\int\limits_{0}^{2\pi}}
\frac{rdrd\theta}{\left(  re^{i\theta}-z_{p}\right)  \left(  re^{-i\theta
}-\overline{z_{p}}\right)  }%
\]%
\[
=\frac{1}{\pi}%
{\displaystyle\int\limits_{0}^{1}}
rdr%
{\displaystyle\int\limits_{0}^{2\pi}}
\frac{d\theta}{\left(  re^{i\theta}-z_{p}\right)  \left(  re^{-i\theta
}-\overline{z_{p}}\right)  }%
\]

Setting $w=e^{i\theta}$ , $dw=iwd\theta,$ $0\leq\theta\leq2\pi$and $\left\vert
w\right\vert =1$%
\[
I_{p}=\frac{1}{\pi}%
{\displaystyle\int\limits_{0}^{1}}
rdr%
{\displaystyle\int\limits_{C_{1}}}
\frac{\frac{dw}{iw}}{\left(  rw-z_{p}\right)  \left(  \dfrac{r}{w}%
-\overline{z_{p}}\right)  }%
\]

i-e%
\[
I_{p}=\frac{1}{\pi i}%
{\displaystyle\int\limits_{0}^{1}}
dr%
{\displaystyle\int\limits_{C_{1}}}
\frac{dw}{\left(  w-\frac{z_{p}}{r}\right)  \left(  r-w\overline{z_{p}%
}\right)  }%
\]

now we state%
\[
\frac{1}{\pi i}%
{\displaystyle\int\limits_{C_{1}}}
\frac{dw}{\left(  w-\frac{z_{p}}{r}\right)  \left(  r-w\overline{z_{p}%
}\right)  }=\frac{2r}{r^{2}-\left\vert z_{p}\right\vert ^{2}}%
\]

hence%
\[
I_{p}=%
{\displaystyle\int\limits_{0}^{1}}
\frac{2r}{r^{2}-\left\vert z_{p}\right\vert ^{2}}dr=\log\left\vert
\frac{1-\left\vert z_{p}\right\vert ^{2}}{\left\vert z_{p}\right\vert ^{2}%
}\right\vert
\]

then
\[%
{\displaystyle\iint\limits_{D}}
\left(
{\displaystyle\sum\limits_{p=1}^{m}}
\left\vert \frac{1}{z-z_{p}}\right\vert ^{2}\right)  dxdy=\pi\log\left\vert
{\displaystyle\prod\limits_{p=1}^{m}}
\frac{1-\left\vert z_{p}\right\vert ^{2}}{\left\vert z_{p}\right\vert ^{2}%
}\right\vert
\]

To prove \eqref{Jpvalue}
\[
J_{p}=\frac{1}{2\pi}%
{\displaystyle\iint\limits_{D}}
\left\vert \frac{z+z_{p}}{z-z_{p}}\right\vert ^{2}dxdy\text{ \ \ \ \ ,
\ \ \ }z=x+iy
\]

suppose \ $z=re^{i\theta}=x+iy$%
\[
dxdy=rdrd\theta,\text{ \ }0\leq r<1,\text{ }0\leq\theta\leq2\pi
\]

then%
\[
J_{p}=\frac{1}{2\pi}%
{\displaystyle\int\limits_{0}^{1}}
{\displaystyle\int\limits_{0}^{2\pi}}
\frac{\left(  re^{i\theta}+z_{p}\right)  \left(  re^{-i\theta}+\overline
{z_{p}}\right)  }{\left(  re^{i\theta}-z_{p}\right)  \left(  re^{-i\theta
}-\overline{z_{p}}\right)  }rdrd\theta
\]%
\[
=\frac{1}{2\pi}%
{\displaystyle\int\limits_{0}^{1}}
rdr%
{\displaystyle\int\limits_{0}^{2\pi}}
\frac{\left(  re^{i\theta}+z_{p}\right)  \left(  re^{-i\theta}+\overline
{z_{p}}\right)  }{\left(  re^{i\theta}-z_{p}\right)  \left(  re^{-i\theta
}-\overline{z_{p}}\right)  }d\theta
\]

Setting $w=e^{i\theta}$ , $dw=iwd\theta,$ $0\leq\theta\leq2\pi$and $\left\vert
w\right\vert =1$%
\[
J_{p}=\frac{1}{2\pi}%
{\displaystyle\int\limits_{0}^{1}}
rdr%
{\displaystyle\int\limits_{C_{1}}}
\frac{\left(  rw+z_{p}\right)  \left(  r+w\overline{z_{p}}\right)  }{\left(
rw-z_{p}\right)  \left(  r-w\overline{z_{p}}\right)  }\frac{dw}{iw}%
\]
i-e%
\[
J_{p}=\frac{1}{2\pi i}%
{\displaystyle\int\limits_{0}^{1}}
dr%
{\displaystyle\int\limits_{C_{1}}}
\frac{\left(  rw+z_{p}\right)  \left(  r+w\overline{z_{p}}\right)  }{\left(
w-\frac{z_{p}}{r}\right)  \left(  r-w\overline{z_{p}}\right)  }\frac{dw}%
{w}=\frac{2r\left(  r^{2}+\left\vert z_{p}\right\vert ^{2}\right)  -r\left(
r^{2}-\left\vert z_{p}\right\vert ^{2}\right)  }{r^{2}-\left\vert
z_{p}\right\vert ^{2}}%
\]

now we state%
\[
J_{p}=%
{\displaystyle\int\limits_{0}^{1}}
r\frac{r^{2}+3\left\vert z_{p}\right\vert ^{2}}{r^{2}-\left\vert
z_{p}\right\vert ^{2}}dr=%
{\displaystyle\int\limits_{0}^{1}}
\left(  r+\frac{4r\left\vert z_{p}\right\vert ^{2}}{r^{2}-\left\vert
z_{p}\right\vert ^{2}}\right)  dr
\]

hence%
\[
J_{p}=\frac{1}{2}+2\left\vert z_{p}\right\vert ^{2}\log\left\vert
\frac{1-\left\vert z_{p}\right\vert ^{2}}{\left\vert z_{p}\right\vert ^{2}%
}\right\vert
\]

then
\[%
{\displaystyle\iint\limits_{D}}
\left(
{\displaystyle\sum\limits_{p=1}^{m}}
\left\vert \frac{z+z_{p}}{z-z_{p}}\right\vert ^{2}\right)  dxdy=m\pi+2\pi
\log\left\vert
{\displaystyle\prod\limits_{p=1}^{m}}
\frac{1-\left\vert z_{p}\right\vert ^{2}}{\left\vert z_{p}\right\vert ^{2}%
}\right\vert ^{2\left\vert z_{p}\right\vert ^{2}}%
\]

and the proposition is proved.
\end{proof}

\section{\bigskip\textbf{Tchebychev }polynomials of the second kind on the
ellipse}

\bigskip$D$ is the ellipse $b^{2}x^{2}+a^{2}y^{2}<a^{2}b^{2};$by the conformal
mapping $z=\cos w$ the cut ellipse is transformed into $R$ the rectangle
$\left(  -ci,-ci+\pi,ci+\pi,ci\right)  $ where $a=\cosh c,b=\sinh c.\left(
c>0\right)  .$We assume that the foci of the ellipse are situated at $z=\pm1$
,that is $a^{2}-b^{2}=1.$

The function $T_{n}\left(  z\right)  =\cos\left(  n\cos^{-1}z\right)
,n=0,1,2,3....$is a polynomial of degree $n$ . $T_{n}\left(  z\right)  $ is
called the Tchebichev polynomial of degree $n$ .

The polynomial of degree $n$
\begin{equation}
U_{n}\left(  z\right)  =\frac{T_{n+1}^{\prime}\left(  z\right)  }{n+1}%
=\frac{\sin\left(  \left(  n+1\right)  \cos^{-1}z\right)  }{\sqrt{1-z^{2}}}
\label{Unz}%
\end{equation}

is called the Tchebichev polynomial of second kind.We shall show that the
polynomials $U_{n}\left(  z\right)  $ are orthogonal to each other with
respect to the above mentioned ellipse.

\begin{theorem}
The orthonormalized polynomials ,$\left\{  P_{n}\left(  z\right)  \right\}
_{n=0,1,2,3.....}$%
\[
P_{n}\left(  z\right)  =2\sqrt{\frac{n+1}{\pi}}\left(  \rho^{n+1}-\rho
^{-n-1}\right)  ^{-\frac{1}{2}}U_{n}\left(  z\right)  \text{ \ \ ,}%
n=0,1,2,3......
\]

satisfies
\begin{equation}%
{\displaystyle\iint\limits_{D}}
P_{n}\left(  z\right)  \overline{P_{m}\left(  z\right)  }dxdy=\delta_{n,m}
\label{Orthon}%
\end{equation}

$\delta_{n,m}$ is the symbol of Kronecker, and%

\begin{equation}
A=%
{\displaystyle\iint\limits_{D}}
\left\vert T_{n}^{\prime}\left(  z\right)  \right\vert ^{2}dxdy=\frac{n\pi}%
{4}\left(  \rho^{n}-\rho^{-n}\right)  \label{Tnarea}%
\end{equation}

\end{theorem}

\begin{proof}
We apply to the ellipse two cuts from $-a$ to $-1$ and from $1$ to $a$
respectively.Obviously ,these cuts do not affect the value of the area
integral
\[
A_{n,m}=%
{\displaystyle\iint\limits_{D}}
U_{n}\left(  z\right)  \overline{U_{m}\left(  z\right)  }dxdy
\]

The conformal mapping $z=\cos w$ the cut ellipse is transformed into the
rectangle $R$ :$\left(  -ci,-ci+\pi,ci+\pi,ci\right)  $,where $a=\cosh c$
,$b=\sinh c$ \ ,$c>0.$Since the Jacobien of the

transformation is%

\begin{equation}
\frac{\partial\left(  x,y\right)  }{\partial\left(  u,v\right)  }=\left\vert
\frac{dz}{dw}\right\vert ^{2}=\left\vert 1-z^{2}\right\vert \label{Jacobien}%
\end{equation}

where%
\[
z=x+iy=\cos w\text{ \ , }w=u+iv
\]

it follows from \eqref{Unz},\eqref{Jacobien} that%
\begin{equation}
A_{n,m}=%
{\displaystyle\iint\limits_{R}}
\sin\left(  \left(  n+1\right)  w\right)  \overline{\sin\left(  \left(
n+1\right)  w\right)  }dudv \label{Anm}%
\end{equation}

Let us evaluate $A_{n-1,m-1}$%
\[
A_{n-1,m-1}=%
{\displaystyle\int\limits_{-c}^{c}}
{\displaystyle\int\limits_{0}^{\pi}}
\sin\left(  n\left(  u+iv\right)  \right)  \overline{\sin\left(  m\left(
u+iv\right)  \right)  }dudv
\]%
\[
=%
{\displaystyle\int\limits_{-c}^{c}}
{\displaystyle\int\limits_{0}^{\pi}}
\sin\left(  nu+inv\right)  \sin\left(  mu-imv\right)  dudv
\]

In fact%
\[
\sin\left(  nu+inv\right)  =\sin nu\cosh nv+i\cos nu\sinh nv
\]

and%
\[
\sin\left(  mu-imv\right)  =\sin mu\cosh mv-i\cos mu\sinh mv
\]

then%
\begin{equation}
A_{n-1,m-1}=%
{\displaystyle\int\limits_{-c}^{c}}
\cosh nv\cosh mvdv%
{\displaystyle\int\limits_{0}^{\pi}}
\sin nu\sin mudu \label{intehr1}%
\end{equation}%
\begin{equation}
+%
{\displaystyle\int\limits_{-c}^{c}}
\sinh nv\sinh mvdv%
{\displaystyle\int\limits_{0}^{\pi}}
\cos nu\cos mudu \label{integr2}%
\end{equation}%
\begin{equation}
+i%
{\displaystyle\int\limits_{-c}^{c}}
\sinh nv\cosh mvdv%
{\displaystyle\int\limits_{0}^{\pi}}
\cos nu\sin mudu \label{integr3}%
\end{equation}%
\begin{equation}
-i%
{\displaystyle\int\limits_{-c}^{c}}
\cosh nv\sinh mvdv%
{\displaystyle\int\limits_{0}^{\pi}}
\sin nu\cos mudu \label{integr4}%
\end{equation}

The last two integrals over $v$ vanish because the integrands are odd.The
first two integrals vanish if $n\neq m.$For $n=m$ we obtain%
\[
A_{n-1,n-1}=%
{\displaystyle\int\limits_{-c}^{c}}
\cosh^{2}nvdv%
{\displaystyle\int\limits_{0}^{\pi}}
\sin^{2}nudu
\]%
\[
+%
{\displaystyle\int\limits_{-c}^{c}}
\sinh^{2}nvdv%
{\displaystyle\int\limits_{0}^{\pi}}
\cos^{2}nudu
\]%
\[
+i%
{\displaystyle\int\limits_{-c}^{c}}
\sinh nv\cosh nvdv%
{\displaystyle\int\limits_{0}^{\pi}}
\cos nu\sin n%
\operatorname{u}%
dudv
\]%
\[
-i%
{\displaystyle\int\limits_{-c}^{c}}
\cosh nv\sinh nvdv%
{\displaystyle\int\limits_{0}^{\pi}}
\sin nu\cos nududv
\]

After some computations
\[
A_{n,n}=\frac{\pi}{2\left(  n+1\right)  }\sinh2\left(  n+1\right)  c
\]

or in view $a+b=\cosh c$ $+\sinh c=e^{c}$
\begin{equation}
A_{n,n}=\frac{\pi}{4\left(  n+1\right)  }\left(  \rho^{n+1}-\rho
^{-n-1}\right)  \text{ \ \ \ \ ,}\left(  a+b\right)  ^{2}=\rho\label{Annn}%
\end{equation}

The orthonormalized polynomials are therefore
\end{proof}

\[
P_{n}\left(  z\right)  =2\sqrt{\frac{n+1}{\pi}}\left(  \rho^{n+1}-\rho
^{-n-1}\right)  ^{-\frac{1}{2}}U_{n}\left(  z\right)
\]

Now by \eqref{Unz},\eqref{Annn},we get
\[%
{\displaystyle\iint\limits_{D}}
\left\vert T_{n}^{\prime}\left(  z\right)  \right\vert ^{2}dxdy=n^{2}%
||U_{n-1}||^{2}=\frac{n\pi}{4}\left(  \rho^{n}-\rho^{-n}\right)
\]

and the theorem is proved.

\section{Interpolation}

We consider not only the variable $z$ but a related complex variable $w$ such
that%
\begin{equation}
z=\frac{1}{2}\left(  w+\frac{1}{w}\right)  \label{Jukovsk}%
\end{equation}
Then, if w moves on the circle $\left\vert w\right\vert =r$ \ for $r>1$
centred at the origin,we have
\[
z=a\cos\theta+ib\sin\theta
\]
where
\begin{equation}
a=\frac{1}{2}\left(  r+\frac{1}{r}\right)  ,\text{ }b=\frac{1}{2}\left(
r-\frac{1}{r}\right)  \label{llipses}%
\end{equation}
Hence z moves on the standard ellipse $E_{r}$ ,with semi axis $a,b.$%
\[
b^{2}x^{2}+a^{2}y^{2}=a^{2}b^{2}%
\]
Namely%
\[
E_{r}=\left\{  z:\left\vert z+\sqrt{z^{2}-1}\right\vert =r\right\}
\]

The roots of $U_{n}\left(  x\right)  $ are all real, distinct, symmetric with
respect to the line $x=0$ and are given by the expression \cite{Dilcher}.%
\begin{equation}
z_{k}=\cos\frac{k\pi}{n+1},\text{ \ \ }k=1,2,...,n. \label{roots}%
\end{equation}

and%
\begin{equation}
U_{n}\left(  z\right)  =2^{n}%
{\displaystyle\prod\limits_{k=1}^{n}}
\left(  z-z_{k}\right)  \label{Unn}%
\end{equation}

We verify that,%
\[
\omega_{n}\left(  z\right)  =%
{\displaystyle\prod\limits_{k=1}^{n}}
\left(  z-z_{k}\right)  =\frac{1}{2^{n}}U_{n}\left(  z\right)
\]%
\begin{equation}
=\frac{1}{2^{n}}\frac{w^{n+1}-w^{-n-1}}{w-w^{-1}} \label{Jukovs}%
\end{equation}

where%
\[
z=\frac{1}{2}\left(  w+\frac{1}{w}\right)
\]

if $\left\vert w\right\vert =r$ ,then,$z\in E_{r}$%
\begin{equation}
\left\vert \omega_{n}\left(  z\right)  \right\vert =\frac{1}{2^{n}}\left\vert
\frac{w^{n+1}-w^{-n-1}}{w-w^{-1}}\right\vert \label{Jukovski}%
\end{equation}

\ Now let $G_{r}$ denote the interior of the ellipse $E_{r}$,and suppose $f$
\ is analytic in $\overline{G_{r}}.$ Let $L_{n-1}\left(  z\right)  $ denote
the interpolating polynomial of degree $n-1$ corresponding to the
interpolation points ,$z_{1},z_{2}....z_{n},$we have%
\begin{equation}
L_{n-1}\left(  z\right)  =\frac{1}{2\pi i}%
{\displaystyle\int\limits_{E_{r}}}
\frac{\omega_{n}\left(  t\right)  -\omega_{n}\left(  z\right)  }{t-z}%
\frac{f\left(  t\right)  }{\omega_{n}\left(  t\right)  }dt\text{
\ \ \ \ \ ,}t\in G_{r} \label{Interpola}%
\end{equation}

\begin{proposition}%
\begin{equation}
L_{n-1}\left(  z\right)  =%
{\displaystyle\sum\limits_{k=1}^{n}}
\frac{\omega_{n}\left(  z\right)  }{z-z_{k}}\frac{f\left(  z_{k}\right)
}{\omega_{n}^{\prime}\left(  z_{k}\right)  } \label{Polyno}%
\end{equation}

and%
\begin{equation}
f\left(  z_{k}\right)  -L_{n-1}\left(  z_{k}\right)  =0\text{ \ \ \ \ \ \ \ ,}%
k=1,2......n \label{Equat}%
\end{equation}

Hence,\ for $-1\leq x\leq1,$%
\begin{equation}
Lim_{n\longrightarrow\infty}f\left(  x\right)  -L_{n-1}\left(  x\right)  =0
\label{Limit}%
\end{equation}

\end{proposition}

\begin{proof}
The image of $\left\{  w:\left\vert w\right\vert =R\right\}  ,R>1$ in the
$w$-complex plane by $z=\frac{1}{2}\left(  w+\dfrac{1}{w}\right)  $ is the
ellipse $E_{R}$ in the $z$-complex plane,with semi axis%
\[
a=\frac{1}{2}\left(  R+R^{-1}\right)  \text{ \ \ \ \ ,and \ \ \ }b=\frac{1}%
{2}\left(  R-R^{-1}\right)  .
\]

If $\omega_{n}\left(  z_{k}\right)  =0,$then%
\[
\underset{z_{k}}{\operatorname{Re}s}\frac{\omega_{n}\left(  t\right)
-\omega_{n}\left(  z\right)  }{\left(  t-z\right)  \omega_{n}\left(  t\right)
}f\left(  t\right)  =-\frac{\omega_{n}\left(  z\right)  }{\left(
z_{k}-z\right)  \omega_{n}^{\prime}\left(  z_{k}\right)  }f\left(
z_{k}\right)
\]

Hence%
\[
\frac{1}{2\pi i}%
{\displaystyle\int\limits_{E_{R}}}
\frac{\omega_{n}\left(  t\right)  -\omega_{n}\left(  z\right)  }{t-z}%
\frac{f\left(  t\right)  }{\omega_{n}\left(  t\right)  }dt=-%
{\displaystyle\sum\limits_{k=1}^{n}}
\frac{f\left(  z_{k}\right)  }{\omega_{n}^{\prime}\left(  z_{k}\right)  }%
\frac{\omega_{n}\left(  z\right)  }{z_{k}-z}%
\]

i-e%
\[
L_{n-1}\left(  z\right)  =%
{\displaystyle\sum\limits_{k=1}^{n}}
\frac{\omega_{n}\left(  z\right)  }{\left(  z-z_{k}\right)  \omega_{n}%
^{\prime}\left(  z_{k}\right)  }f\left(  z_{k}\right)
\]

Since%
\begin{equation}
f\left(  z\right)  =\frac{1}{2\pi i}%
{\displaystyle\int\limits_{E_{R}}}
\frac{f\left(  t\right)  }{t-z}dt\text{ } \label{Cauch}%
\end{equation}

By \eqref{Cauch}and\eqref{Interpola},we get%
\begin{equation}
f\left(  z\right)  -L_{n-1}\left(  z\right)  =\frac{1}{2\pi i}%
{\displaystyle\int\limits_{E_{R}}}
\frac{\omega_{n}\left(  z\right)  }{\omega_{n}\left(  t\right)  }%
\frac{f\left(  t\right)  }{t-z}dt\text{ \ \ \ \ \ ,}t\in G_{R}
\label{Difftinterpo}%
\end{equation}

it follows that%
\[
f\left(  z_{k}\right)  -L_{n-1}\left(  z_{k}\right)  =0\text{ \ \ \ \ \ \ \ ,}%
k=1,2......n
\]

Further,if $\left\vert w\right\vert =R,$ \eqref{Jukovski} implies that
\begin{equation}
\left\vert \omega_{n}\left(  z\right)  \right\vert =\frac{R^{n}}{2^{n}%
}\left\vert \frac{1-w^{-2n-2}}{1-w^{-2}}\right\vert \text{ \ \ \ , \ \ \ if
}z\in E_{R} \label{Image}%
\end{equation}

\ For $-1\leq x\leq1,$relations \eqref{Image}and \eqref{Difftinterpo}
immediately imply the estimate%
\[
f\left(  x\right)  -L_{n-1}\left(  x\right)  =\frac{1}{2\pi i}%
{\displaystyle\int\limits_{E_{R}}}
\frac{\omega_{n}\left(  x\right)  }{\omega_{n}\left(  t\right)  }%
\frac{f\left(  t\right)  }{t-x}dt\text{ \ \ }%
\]%
\[
=O\left(  1\right)  .\frac{1}{2^{n}}.\frac{2^{n}}{R^{n}}=O\left(  \frac
{1}{R^{n}}\right)
\]

As $n\longrightarrow\infty.$ The larger the domain surrounding $\left[
-1\text{ \ \ }1\right]  $ is in which $f$ \ is analytic , the faster the
convergence $L_{n}\Longrightarrow f$ \ takes place in $\left[  -1\text{
\ \ }1\right]  .$
\end{proof}

\end{document}